\documentclass{amsart}

\usepackage[utf8]{inputenc}
\usepackage[english]{babel}
\usepackage[babel=true,kerning,spacing]{microtype}
\usepackage{amssymb,amsfonts,amsmath,amsthm}
\usepackage{url}
\usepackage{graphicx}
\renewcommand{\le}{\leqslant}
\renewcommand{\ge}{\geqslant}
\renewcommand{\leq}{\leqslant}
\renewcommand{\geq}{\geqslant}
\usepackage{tikz}

\newtheorem*{conjecture}{Conjecture}
\newtheorem*{theorem}{Theorem}
\newtheorem*{proposition}{Proposition}
\newtheorem*{corollary}{Corollary}

\newcommand{\Prob}{\operatorname{Prob}}

\newcommand{\UU}{\mathbb U}

\newcommand{\ZZ}{\mathbb Z}

\newcommand{\FF}{\mathbb F}
\renewcommand{\AA}{\mathcal A}
\newcommand{\BB}{\mathcal B}
\newcommand{\CC}{\mathcal C}
\newcommand{\PS}{\mathcal P}
\newcommand{\E}{E}
\newcommand{\N}{\operatorname{N}}

\newcommand{\cl}[1]{\operatorname{cl}(#1)}

\newcommand{\ord}{\mathcal{O}}
\newcommand{\End}[1]{\operatorname{End}(#1)}

\newcommand{\Ell}[1]{\operatorname{Ell}(#1)}

\newcommand{\pimap}{\pi}
\newcommand{\hmap}{\eta}
\newcommand{\tphi}{\varphi}
\newcommand{\gltwo}[1]{\operatorname{GL}(2,#1)}
\newcommand{\rhotot}{\rho_\text{tot}}
\newcommand{\Exp}[1]{\operatorname{\mathbf{E}}[#1]}
\title[Finding short product representations in finite groups]{A low-memory
algorithm for finding short product representations in finite groups}
\author{Gaetan Bisson and Andrew V. Sutherland}
\date{}

\begin{document}

\begin{abstract}
We describe a space-efficient algorithm for solving a generalization of the
subset sum problem in a finite group $G$, using a Pollard-$\rho$ approach.
Given an element $z$ and a sequence of elements $S$, our algorithm attempts to
find a subsequence of $S$ whose product in $G$ is equal to~$z$. For a random
sequence $S$ of length $d\log_2 n$, where $n=\#G$ and $d \ge 2$ is a constant,
we find that its expected running time is $O(\sqrt{n}\log n)$ group operations
(we give a rigorous proof for $d > 4$), and it only needs to store $O(1)$
group elements.  We consider applications to class groups of imaginary quadratic
fields, and to finding isogenies between elliptic curves over a finite field.
\end{abstract}
\maketitle

%%%%%%%%%%%%%%%%%%%%%%%%%%%%%%%%%%%%%%%%%%%%%%%%%%%%%%%%%%%%%%%%
\section{Introduction}

Let $S$ be a sequence of elements in a finite group $G$ of order $n$, written multiplicatively.
We say that $S$
\emph{represents} $G$ if every element of $G$ can be expressed as the (ordered)
product of a subsequence of~$S$. Ideally, we want~$S$ to be short, say
$k=d\log_2 n$ for some constant $d$ known as the \emph{density} of $S$.

In order for $S$ to represent $G$, we clearly require $d\ge 1$, and for
sufficiently large~$n$, any $d>1$ suffices. More precisely, Babai and
Erd\H{o}s~\cite{babai-erdos} show that for all
\[
k \ge \log_2 n + \log_2 \log n + 2
\]
there exists a sequence $S$ of length $k$ that represents~$G$. Their proof is
non-constructive, but, in the case that $G$ is abelian, Erd\H{o}s and
R\'{e}nyi~\cite{erdos-renyi} show that a randomly chosen sequence of length
\[
k = \log_2 n + \log_2 \log n + \omega_n
\]
represents $G$ with probability approaching $1$ as $n\to\infty$, provided that
$\omega_n\to\infty$. The randomness assumption is necessary, since it takes
much larger values of $k$ to ensure that \emph{every} sequence of length $k$
represents $G$, see~\cite{eggleton-erdos,white}.

In related work, Impagliazzo and Naor prove that for a random sequence~$S$ of
density $d>1$, the distribution of subsequence products almost surely
converges to the uniform distribution on $G$ as $n$ goes to
infinity~\cite[Proposition~4.1]{impagliazzo-naor}.  This result allows us to
bound the complexity of our algorithm for almost all~$S$
with $d > 4$.

\medskip

Given a sequence $S$ that represents $G$ (or a large subset of $G$), we
wish to find an explicit representation of a given group element $z$ as the
product of a subsequence of $S$; we call this a \emph{short product
representation} of~$z$. In the special case that $G$ is abelian and the
elements of $S$ are distinct, this is the \emph{subset sum problem} in a finite
group. Variations of this problem and its decision version have long been of
interest to many fields: complexity theory~\cite{karp},
cryptography~\cite{merkle-hellman}, additive number theory~\cite{babai-erdos},
Cayley graph theory~\cite{alon-milman}, and information
theory~\cite{alon-barak-manber}, to name just a few.

As a computational framework, we work with a generic group $G$ whose elements
are uniquely identified, and assume that all group operations are performed by
a black box that can also provide random group elements;
see~\cite[Chapter~1]{sutherland:thesis} for a formal model. Time complexity is
measured by counting group operations (calls to the black box), and for space
complexity we count the number of group elements that are simultaneously
stored. In most practical applications, these metrics are within a
polylogarithmic factor of the usual bit complexity.

Working in this model ensures that our algorithms apply to any
finite group for which a suitable black box can be constructed. It also means
that finding short product representations is provably hard.
Indeed, the discrete logarithm problem in a cyclic group of prime order has a
lower bound of $\Omega(\sqrt{n})$ in the generic group model~\cite{shoup},
and is easily reduced to finding short product representations.

In the particular group $G=\ZZ/n\ZZ$, we note that finding short product
representations is easier for non-generic algorithms: the problem can be
lifted to $k$ subset sum problems in $\ZZ$, which for suitable inputs can be
solved with a time and space complexity of $O(n^{0.3113})$
via~\cite{howgravegraham-joux}, beating the $\Omega(\sqrt{n})$ generic lower bound
noted above. This is not so surprising, since working with integers is often
easier than working in generic groups; for instance, the discrete logarithm
problem in $\ZZ$ corresponds to integer division and can be solved in
quasi-linear time.

\medskip

A standard technique for solving subset sum problems in generic groups uses a
baby-step giant-step approach, which can also be used to find short product
representations (Section~\ref{sec:BSGS}). This typically involves $O(2^{k/2})$
group operations and storage for $O(2^{k/2})$ group elements. The space bound
can be improved to $O(2^{k/4})$ via a method of Schroeppel and
Shamir~\cite{schroeppel-shamir}.

Here, we give a Pollard-$\rho$ type algorithm~\cite{pollard} for finding short
product representations in a finite group (Section~\ref{sec:pollard}).
It only needs to store $O(1)$ group elements, and, assuming $S$ is a random
sequence of density $d>4$, we prove that its expected running time is
$O(\sqrt{n}\log{n})$ group operations; alternatively, by dedicating $O(n^\epsilon)$
space to precomputations, the time complexity can be reduced to $O(\sqrt{n})$
(Section~\ref{sec:analysis}).

We also consider two applications: representing elements of the class group of
an imaginary quadratic number field as short products of prime ideals with
small norm (Section~\ref{sec:relations}), and finding an isogeny between two
elliptic curves defined over a finite field (Section~\ref{sec:isogenies}). For
the latter, our method combines the advantages of~\cite{galbraith}
and~\cite{galbraith-hess-smart} in that it requires little memory and finds an
isogeny that can subsequently be evaluated in polynomial time.

In practice, our algorithm performs well so long as $d \ge 2$, and its low
space complexity allows it to feasibly handle much larger problem instances
than other generic methods (Section~\ref{sec:comput}).

%%%%%%%%%%%%%%%%%%%%%%%%%%%%%%%%%%%%%%%%%%%%%%%%%%%%%%%%%%%%%%%%
\section{Algorithms}

Let $S$ be a sequence of length $k$ in a finite group $G$ of order $n$, let $z$
be an element of $G$, and let $\PS(S)$ denote the set of all subsequences
of~$S$. Our goal is to find a preimage of $z$ under the product map
$\pimap:\PS(S)\to G$ that sends a subsequence of~$S$ to the (ordered) product of
its elements.

\subsection{Baby-step giant-step}\label{sec:BSGS}

Let us first recall the baby-step giant-step method. We may
express~$S=AB$ as the concatenation of two subsequences of roughly equal length.
For any sequence $y=(y_1,\ldots,y_m)$, let $\mu(y) = (y_m^{-1},\ldots,y_1^{-1})$,
so that $\pimap(y)$ and $\pimap(\mu(y))$ are inverses in~$G$. We then search for
$x\in\PS(A)$ (a baby step) and $y\in\PS(B)$ (a giant step)
which ``collide'' in the sense that $\pimap(x) = \pimap(z\mu(y))$,
where $z\mu(y)$ denotes the sequence $(z,y_m^{-1},\ldots,y_1^{-1})$.

\smallskip
\begin{quote}
\textsc{Baby-step giant-step Algorithm}\\
\textbf{\scshape Input:} A finite sequence $S$ in a group $G$ and a target $z\in\pimap(\PS(S))$.\\
\textbf{\scshape Output:} A subsequence of $S$ whose product is~$z$.\\
\begin{tabular}{rl}
   1. & Express $S$ in the form $S=AB$ with $\#A\approx \#B$.
\\ 2. & For each $x\in\PS(A)$, store $(\pimap(x),x)$ in a table indexed by $\pimap(x)$.
\\ 3. & For each $y\in\PS(B)$:
\\ 4. & \qquad Lookup $\pimap(z\mu(y))$ in the table computed in Step~2.
\\ 5. & \qquad If $\pimap(z\mu(y))=\pimap(x)$ is found then output $xy$, otherwise continue.
\end{tabular}
\end{quote}
\smallskip

The table constructed in Step~2 is typically implemented as a hash table, so
that the cost of the lookup in Step~4 is negligible.  Elements of $\PS(A)$ and
$\PS(B)$ may be compactly represented by bit-strings of length $\lceil
k/2\rceil = O(\log n)$, which is approximately the size of a single group
element. If these bit-strings are enumerated in a suitable order, each step can
be derived from the previous step using $O(1)$ group operations\footnote{With a Gray code,
exactly one group operation is used per step,
see~\cite{knuth-art4f2}.}. The algorithm then performs a total of $O(2^{k/2})$ group
operations and has a space complexity of $O(2^{k/2})$ group elements. One can
make a time-space trade off by varying the relative sizes of $A$ and $B$.

This algorithm has the virtue of determinism, but its
complexity $O(n^{d/2})$ is exponential in the density $d$ (as well as $\log n$). For $d > 1$, a
randomized approach works better:  select $\sqrt{n}$ baby steps $x\in\PS(A)$
at random, then select random giant steps $y\in\PS(B)$ until a collision
$\pimap(z\mu(y))=\pimap(x)$ is found.  Assuming that
$\pimap(x)$ and $\pimap(z\mu(y))$ are uniformly distributed in $G$,
we expect to use $\sqrt{n}$ giant steps.  To reduce the cost of each step,
one may partition $A$ and $B$ each into approximately $d$
subsequences $A_i$ and $B_i$ and precompute $\pimap(x)$ for all $x\in\PS(A_i)$,
and $\pimap(\mu(y))$ for all $y\in\PS(B_i)$. This yields an expected running
time of $O(\sqrt{n})$ group operations, using storage for $O(\sqrt{n})$ group
elements, for any fixed $d$.

\subsection{A low-memory algorithm}\label{sec:pollard}

In order to use the Pollard-$\rho$ technique, we need a pseudo-random function
$\phi$ on the disjoint union $\CC=\AA\sqcup\BB$, where $\AA=\PS(A)$ and $\BB$ is the set
$\{z\mu(y):y\in\PS(B)\}$.  This map $\phi$ is required to preserve collisions, meaning that
$\pimap(x)=\pimap(y)$ implies $\pimap(\phi(x))=\pimap(\phi(y))$. Given a hash function $\hmap:G\to\CC$, we may construct such a map as $\phi=\hmap\circ\pimap$. Under suitable assumptions (see Section~\ref{sec:analysis}), the Pollard-$\rho$ method can then be applied.

\smallskip
\begin{quote}
\textsc{Pollard-$\rho$ Algorithm}\\
\textbf{\scshape Input:} A finite sequence $S$ in a group $G$ and a target $z\in \pimap(\PS(S))$.\\
\textbf{\scshape Output:} A subsequence of $S$ whose product is~$z$.\\
\begin{tabular}{rl}
   1. & Pick a random element $w\in \CC$ and a hash function $\hmap:G\to \CC$.
\\ 2. & Find the least $i > 0$ and $j \ge 0$ such that $\phi^{(i+j)}(w)=\phi^{(j)}(w)$.
\\ 3. & If $j=0$ then return to Step~1.
\\ 4. & Let $s=\phi^{(i+j-1)}(w)$ and let $t=\phi^{(j-1)}(w)$.
\\ 5. & If $\pimap(s) \ne \pimap(t)$ then return to Step~1.
\\ 6. & If $s\in\AA$ and $t=z\mu(y)\in\BB$ then output $sy$ and terminate.
\\ 7. & If $t\in\AA$ and $s=z\mu(y)\in\BB$ then output $ty$ and terminate.
\\ 8. & Return to Step~1.
\end{tabular}
\end{quote}
\smallskip

Step~2 can be implemented with Floyd's
algorithm~\cite[Exercise~3.1.6]{knuth-art2} using storage for just two elements of $\CC$,
which fits in the memory space of $O(1)$ group elements.
More sophisticated collision-detection techniques can reduce the number
of evaluations of $\phi$ while still storing $O(1)$ elements,
see~\cite{brent,sedgewick,teske}.  We prefer the method of \emph{distinguished points},
which facilitates a parallel implementation~\cite{vanoorschot-wiener}.

\subsection{Toy example}

Let $G=(\ZZ/n\ZZ,+)$ and define $S$ as the concatenation of the
sequences $A=(3^i)$ and $B=(5^i)$ for $i\in\{1,\ldots,k/2\}$. We put $n=127$ and
$k=12$, implying $d\approx 1.7$. With $\CC=\AA\sqcup\BB$ as above, we define
$\hmap:G\to \CC$ via
\[
x\longmapsto\left\{\begin{array}{cl}
(A_i)_{\{i:b_i=1\}}&\text{when }b_0=1\\
z\mu\left((B_i)_{\{i:b_i=1\}}\right)&\text{when }b_0=0
\end{array}\right.
\]
where $\sum_{i=0}^{k/2} b_i2^i$ is the binary representation of $96x \bmod n$.

Starting from $w=(2,-5^6,-5^3,-5^2,-5)$, the algorithm finds $i=4$ and $j=6$:
\smallskip

\begin{center}
\tiny
\begin{tikzpicture}
  \foreach \i/\s/\x/\y in {
    0/{2,-5^6,-5^3,-5^2,-5}/0/0.85,
    1/{3^3,3^5}/0/1.7,
    2/{2,-5^5,-5^4}/2/1.7,
    3/{2,-5^6,-5^5,-5^4,-5^2,-5}/5/1.7,
    4/{3^2,3^4}/7.8/1.7,
    5/{2,-5^5}/9.5/1.7,
    6/{3,3^2,3^5}/10.5/0.85,
    7/{2,-5^2,-5}/10/0,
    8/{2,-5^6,-5^4,-5^2,-5}/7.3/0,
    9/{3,3^2,3^3,3^5}/6.8/0.85}
  \node (G-\i) at (\x,\y) {$\left(\s\right)$};
  \foreach \from/\to in {0/1,1/2,2/3,3/4,4/5,5/6,6/7,7/8,8/9,9/4}
    \draw[black,->] (G-\from) -- (G-\to);
\end{tikzpicture}
\normalsize
\end{center}
\smallskip

\noindent
The two preimages of $(3^2,3^4)$ yield the short product representation
\[
2\equiv 3+3^2+3^3+3^5+5+5^2+5^4+5^5+5^6\bmod 127.
\]

%%%%%%%%%%%%%%%%%%%%%%%%%%%%%%%%%%%%%%%%%%%%%%%%%%%%%%%%%%%%%%%%
\section{Analysis}\label{sec:analysis}

The Pollard-$\rho$ approach is motivated by the following observation: if
$\phi:X\to X$ is a random function on a set $X$ of cardinality $n$, then the
expected size of the orbit of any $x\in X$ under the action of $\phi$ is
$\sqrt{\pi n/2}$ (see~\cite{sobol-random-sequences} for a rigorous proof). In
our setting, $X$ is the set $\CC$ and $\phi=\hmap\circ\pimap$. Alternatively, since
$\phi$ preserves collisions, we may regard $X$ as the set $\pimap(\CC)\subset G$ and
use $\tphi=\pimap\circ \hmap$.  We shall take the latter view, since it simplifies our
analysis.

Typically the function $\tphi$ is not truly random,
but under a suitable set of assumptions it may behave so. To
rigorously analyze the complexity of our algorithm, we fix a real number $d>4$
and assume that:
\begin{enumerate}
\item the hash function $\hmap:G\to \CC$ is a random oracle;
\item $S$ is a random sequence of density~$d$.
\end{enumerate}

For any finite set $U$, let $\UU_U$ denote the uniform distribution on $U$,
which assigns to each subset $X$ of $U$ the value $\#X/\#U$. For any function
$f:U\to V$, let $f_*\UU_U$ denote the \emph{pushforward distribution} by $f$ of
$\UU_U$, which assigns to each subset~$Y$ of~$V$ the value
\[
f_*\UU_U(Y) = \frac{\#\{u\in U: f(u)\in Y\}}{\#U}.
\]

Assumption~(2) implies that $A$ and $B$ are both random sequences with density
greater than~$2$. By~\cite[Proposition~4.1]{impagliazzo-naor}, this implies
that
\[
\Prob_A\left[
\left\|\pimap_*\UU_{\AA}-\UU_{G}\right\|\geq n^{-c}
\right]\leq n^{-c},
\]
where $c=(d-2)/4 > 1/2$, and the \emph{variation distance} $\|\sigma-\tau\|$
between two distributions $\sigma$ and $\tau$ on $G$ is defined as the maximum
value of $|\sigma(H)-\tau(H)|$ over all subsets $H$ of~$G$. Similarly, we have
\[
\Prob_B\left[
\left\|\pimap_*\UU_{\BB}-\UU_{G}\right\|\geq n^{-c}
\right]\leq n^{-c}.
\]

{}From now on we assume that $S$ is fixed and that $\pimap_*\UU_C$ is within
variation distance $2n^{-c}$ of the uniform distribution on $G$; by the
argument above, this happens with probability at least~$1-2n^{-c}$. Recall that
a \emph{random oracle} $\hmap:G\to \CC$ is a random function drawn uniformly from
$\CC^G$, that is, each value $\hmap(x)$ is drawn uniformly and independently from
$\CC$.  Thus, for any $g\in G$, the distribution of $\pimap(\hmap(g))$ is $\pimap_*\UU_C$.
It is then easy to verify that
\[
\left\|(\hmap\mapsto \pimap\circ \hmap)_*\UU_{\CC^G}-\UU_{G^G}\right\|\leq 2n^{-c}.
\]
In other words, for a random oracle $\hmap$, the function $\tphi=\pimap\circ \hmap$ is very
close to being a random oracle (from $G$ to $G$) itself.

Since $c>1/2$, we obtain, as in~\cite{pollard}, an $O(\sqrt{n})$ bound on the
expectation of the least positive integer $i+j$ for which
$\tphi^{(i+j)}(g)=\tphi^{(j)}(g)$, for any $g=\pi(w)\in G$.  For $d > 2$, the probability
that $\pimap(s)\ne \pimap(t)$ in Step~5 is $o(1)$, since $\CC$ is then
larger than~$G$ and collisions in the map $\tphi$ (and $\phi$) are more likely to be caused by collisions in
$\pimap$ than collisions in $\hmap$.  Having reached Step~6, 
we obtain a short product representation of~$z$ with probability $1/2$, since by
results of~\cite{impagliazzo-naor} the value of $\pimap(x)$ is independent of whether
$x\in \AA$ or $x\in \BB$. The expected running time
is thus $O(k\sqrt{n})=O(\sqrt{n}\log n)$ group operations, and, as noted in
Section~\ref{sec:pollard}, the space complexity is $O(1)$ group elements. We summarize our
analysis with the following proposition.

\begin{proposition}\label{prop:main}
Let $S$ be a random sequence of constant density $d > 4$ and let $\hmap:G\to \CC$ be
a random oracle. Then our Pollard-$\rho$ algorithm
uses $O(\sqrt{n}\log n)$ expected group operations and storage for $O(1)$ group
elements.
\end{proposition}

As in Section~\ref{sec:BSGS}, to speed up the evaluation of the product map
$\pimap$, one may partition $A$ and $B$ into subsequences $A_i$ and $B_i$ of
length $m$ and precompute $\pimap(\PS(A_i))$ and $\pimap(\mu(\PS(B_i))$. This
requires storage for $O(k2^m/m)$ group elements and speeds up subsequent
evaluations of $\pimap$ by a factor of $m$. If we let $m=\epsilon\log_2 n$, for
any $\epsilon>0$, we obtain the following corollary.

\begin{corollary}
Under the hypotheses of the proposition above, our Pollard-$\rho$
algorithm can be implemented to run in expected
time $O(\sqrt{n})$ using $O(n^\epsilon)$ space.
\end{corollary}

In our analysis above, we use a random $S$ random with $d > 4$ to prove that products of
random elements of $\AA$ and $\BB$ are quasi-uniformly distributed in $G$. If
we directly assume that both $\pimap_*\UU_\AA$ and $\pimap_*\UU_\BB$ are
quasi-uniformly distributed, our analysis applies to all $d\ge 2$, and in
practice we find this to be the case. However, we note that this does not apply
to $d<2$, for which we expect a running time of $O(n^{(4-d)/4}\log n)$, as
discussed in Section~\ref{sec:comput}.

%%%%%%%%%%%%%%%%%%%%%%%%%%%%%%%%%%%%%%%%%%%%%%%%%%%%%%%%%%%%%%%%
\section{Applications}

As a first application, let us consider the case where $G$ is the ideal class
group of an order $\ord$ in an imaginary quadratic field.  We may assume
\[
\ord=\ZZ+\frac{D+\sqrt{D}}{2}\ZZ,
\]
where the \emph{discriminant} $D$ is a negative integer congruent to $0$ or $1$
modulo~$4$. Modulo principal ideals, the invertible ideals of $\ord$ form a
finite abelian group $\cl\ord$ of cardinality $h$. The \emph{class number} $h$
varies with $D$, but is on average proportional to $\sqrt{|D|}$ (more
precisely, $\log h \sim \frac{1}{2}\log|D|$ as $D\to -\infty$, by Siegel's
theorem~\cite{siegel}).
Computationally, invertible $\ord$-ideals can be represented as binary
quadratic forms, allowing group operations in $\cl\ord$ to be computed in
time $O(\log^{1+\epsilon}|D|)$, via~\cite{schonhage-fastforms}.

\subsection{Prime ideals}\label{sec:Sk}

Let $\ell_i$ denote the $i$\textsuperscript{th} largest prime number for which
there exists an invertible $\ord$-ideal of norm $\ell_i$ and let $\alpha_i$ denote
the unique such ideal that has nonnegative trace.  For each positive
integer $k$, let $S_k$ denote the sequence of (not necessarily distinct) ideal
classes
\[
S_k = ([\alpha_1],[\alpha_2],\ldots,[\alpha_k]).
\]
For algorithms that work with ideal class groups, $S_k$ is commonly used as a
set of generators for $\cl\ord$, and in practice $k$ can be made quite small,
conjecturally $O(\log h)$.  Proving such a claim is believed to be very
difficult, but under the generalized Riemann hypothesis (GRH), Bach obtains the
following result~\cite{bach-erh}.

\begin{theorem}[Bach]
Assume the GRH. If $D$ is a fundamental\footnote{Meaning that either $D$ is
square-free, or $D/4$ is an integer that is square-free modulo~$4$.} discriminant
and $\ell_{k+1} > 6\log^2|D|$, then the set $S_k$ generates $\cl\ord$.
\end{theorem}

Unfortunately, this says nothing about short product representations in
$\cl\ord$. Recently, a special case of~\cite[Corollary~1.3]{expander-grh}
was considered in~\cite[Theorem~2.1]{quantum-iso} which still assumes the GRH
but is more suited to our short product representation setting.
Nevertheless, for our purpose here, we make the following stronger conjecture.

\begin{conjecture}
For every $d_0 >1$ there exist constants $c > 0$ and $D_0 < 0$ such that if $D
\leq D_0$ and $S_k$ has density $d \geq d_0$ then
\begin{enumerate}
\item $\pimap(\PS(S_k))=G$, that is, $S_k$ represents $G$;
\item $\left\|\pimap_*\UU_{\PS(S_k)}-\UU_G\right\|<h^{-c}$;
\end{enumerate}
where $G$ is the ideal class group $\cl\ord$ and $h$ is its cardinality.
\end{conjecture}

In essence, these are heuristic analogs to the results of Erd\H{o}s and
R\'{e}nyi, and of Impagliazzo and Naor, respectively, suggesting that the
distribution of the classes~$[\alpha_i]$ resembles that of random elements
uniformly drawn from $\cl\ord$. Note that (1), although seemingly weaker, is
only implied by (2) when $c>1$.

\begin{figure}
\begin{center}
\includegraphics[width=\textwidth]{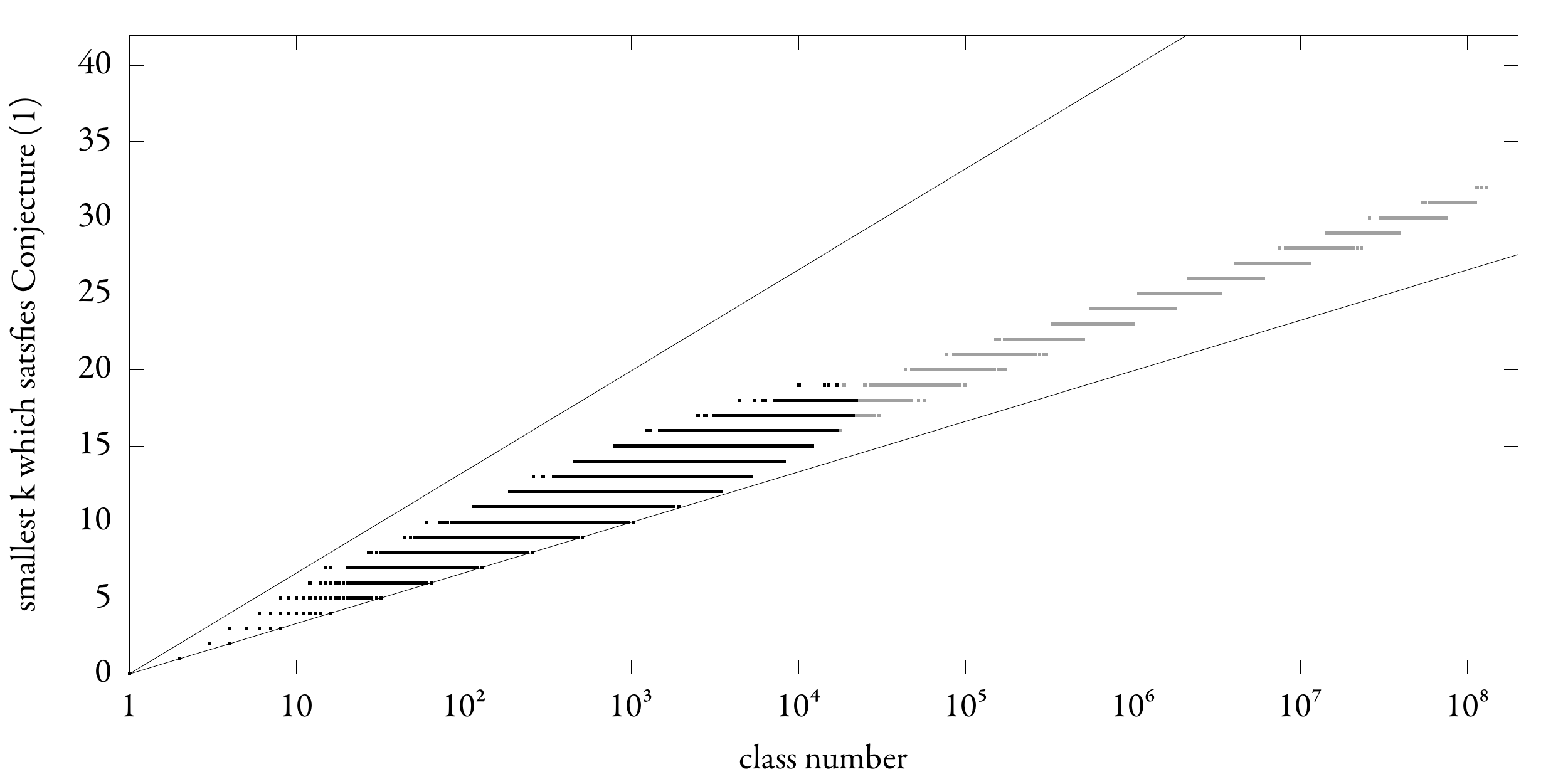}
\end{center}
\caption{\small Dots plot the minimal $k$ such that $S_k$ satisfies conjecture~(1);
gray dots for all discriminants $D\geq -10^8$ and black dots for ten
thousand $D$ drawn at random according to a logarithmic distribution. The lines
represent $k=d\log_2 h$ for $d=1,2$.}
\label{fig:hyp-rnd}
\end{figure}

Empirically, (1) is easily checked: for $d_0=2$ we have verified it using
$D_0=-3$ for every imaginary quadratic order with discriminant $D\geq -10^8$, and for $10^4$
randomly chosen orders with $D$ logarithmically distributed over the
interval $[-10^{16},-10^{8}]$ (see Figure~\ref{fig:hyp-rnd}).  Although harder to
test, (2) is more natural in our context, and practical computations support it
as well. Even though we see no way to prove this conjecture, we assume its
veracity as a useful heuristic.

\subsection{Short relations}\label{sec:relations}

In~\cite{hafner-mccurley}, Hafner and McCurley give a subexponential algorithm
to find representatives of the form $\prod\alpha_i^{e_i}$ for arbitrary ideal
classes of imaginary quadratic orders; the ideals $\alpha_i$ have subexponential norms,
but the exponents $e_i$ can be as large as the class number~$h$.

Asking for small exponents $e_i\in\{0,1\}$ means, in our terminology, writing
elements $z\in G$ as short product representations on $S_k=([\alpha_i])$. Under
the conjecture above, this can be achieved by our low-memory
algorithm in $O(|D|^{1/4+\epsilon})$ expected time, using $k=O(\log h)$ ideals
$\alpha_i$.

We can even combine these approaches.
If the target element $z$ is represented by an ideal of small
norm, say $z=[\alpha_{k+1}]$, we get what we call a \emph{short relation} for~$\cl\ord$.
Conjecture (1) implies not only that the map that sends
each vector $(e_1,\ldots,e_{k+1})\in\ZZ^{k+1}$ to the class of the ideal $\prod\alpha_i^{e_i}$
is surjective, but also that there exists a set
of short relations generating its kernel lattice $\Lambda$. This gives
a much better upper bound on the diameter of $\Lambda$ than was used by Hafner
and McCurley, and their algorithm can be adapted to make use of this new bound
and find, in subexponential time, representatives $\prod\alpha_i^{e_i}$ with
ideals $\alpha_i$ of subexponential norm and exponents $e_i$ bounded by $O(\log|D|)$.
See~\cite{bisson-grh} for details, or~\cite{quantum-iso} for an equivalent construction.

\subsection{Short isogenies}\label{sec:isogenies}

Now let us consider the problem of finding an isogeny between two ordinary elliptic
curves $\E_1$ and $\E_2$ defined over a finite field $\FF_q$.  This problem is of
particular interest to cryptography because the discrete logarithm problem can then
be transported from $\E_1$ to $\E_2$.  An isogeny between curves $\E_1$ and $\E_2$
exists precisely when $\E_1$ and $\E_2$ lie in the same \emph{isogeny class}.  By
a theorem of Tate, this occurs if and only if $\#\E_1(\FF_q)=\#\E_2(\FF_q)$, which
can be determined in polynomial time using Schoof's algorithm~\cite{schoof-pointcounting}.

The isogeny class of $\E_1$ and $\E_2$ can be partitioned according to
the endomorphism rings of the curves it contains, each of which is isomorphic
to an order $\ord$ in an imaginary quadratic number field.  Identifying isomorphic
curves with their $j$-invariant, for each order $\ord$ we define
\[
\Ell\ord=\left\{j(\E) : \End\E\cong\ord\right\},
\]
where $\E$ denotes an elliptic curve defined over $\FF_q$.
The set $\Ell\ord$ to which a given curve belongs can be determined in
subexponential time, under heuristic assumptions~\cite{bisson-sutherland}.
An isogeny from $\E_1$ to $\E_2$ can always be decomposed into two isogenies,
one that is essentially determined by $\End{\E_1}$ and $\End{\E_2}$ (and can
be made completely explicit but may be difficult to compute), and another
connecting curves that lie in the same set $\Ell\ord$.  We shall thus restrict
ourselves to the problem of finding an isogeny between two elements of $\Ell\ord$.

The theory of complex multiplication states that $\Ell\ord$ is a principal
homogeneous space (a \emph{torsor}) for the class group $\cl\ord$: each
ideal $\alpha$ acts on $\Ell\ord$ via an isogeny of degree $\N(\alpha)$, and
this action factors through the class group.  We may then identify each ideal class
$[\alpha]$ with the image $[\alpha]j(\E_i)$
of its action on $j(\E_i)$.  This allows us to effectively work in the group
$\cl\ord$ when computing isogenies from $\E_i$.

Galbraith addressed the search for an isogeny $\E_1\to\E_2$ using a baby-step
giant-step approach in~\cite{galbraith}; a low-memory variant was later given
in~\cite{galbraith-hess-smart} which produces an exponentially long chain of
low-degree isogenies. From that, a linearly long chain of isogenies of
subexponential degree may be derived by smoothing the corresponding ideal in $\cl\ord$ using variants of the method of Hafner and
McCurley (for instance, those mentioned in Section~\ref{sec:relations}); alternatively, our low-memory
algorithm can be used to derive a chain of low-degree isogenies with length
linear in $\log|D|$ (assuming our conjecture), and we believe this is the most
practical approach. However, let us describe how our method applies naturally
to the torsor $\Ell\ord$, and directly finds a short chain of low-degree
isogenies from $\E_1$ to $\E_2$ using very little memory.

\medskip

Let $S_k=AB$ be such that conjecture (1) holds, where $A$ and $B$
are roughly equal in size, and define $\CC=\AA\sqcup\BB$ where $\AA=\PS(A)$ and
$\BB=\mu(\PS(B))$.  We view each element of $\AA$ as a short chain of isogenies of small
prime degree $\ell_i=\N(\alpha_i)$
that originates at $\E_1$; similarly, we view elements of $\BB$ as chains
of isogenies originating at $\E_2$.  Now let $\pimap:\CC\to \Ell\ord$
be the map that sends $x\in\AA$ (resp. $x\in\BB$) to the element of $\Ell\ord$ that is
the codomain of the isogeny chain defined by $x$ and originating at $\E_1$ (resp. $\E_2$).
It suffices to find a collision between an element of $\AA$ and an element
of $\BB$ under the map $\pimap$: this yields an isogeny chain from~$\E_1$ and an isogeny
chain from~$\E_2$ that have the same codomain.  Composing the first with the dual of the
second gives an isogeny from $\E_1$ to $\E_2$.

The iteration function $\phi$ on $\CC$ can now be defined as the composition
$\hmap\circ\pimap$ where~$\hmap$ is a map from $\Ell\ord$ to $\CC$ that behaves like
a random oracle. Using this formalism, our Pollard-$\rho$ algorithm can be
applied directly, and under the conjecture it finds an isogeny in time
$O(h^{1/2+\epsilon})$. In terms of space, it only needs to store $O(1)$
elements of $\cl\ord$ and $\Ell\ord$, which is $O(\log q)$ bits.  However, in
order to compute isogenies, modular polynomials $\Phi_\ell(X,Y)$ might be used,
each of which requires $O(\ell^3\log\ell)$ bits. If we heuristically assume
that $\ell_k = O(k\log k) = O(\log h\log\log h)$, the overall space complexity
is then bounded by $O(\log^{3+\epsilon} h) = O(\log^{3+\epsilon} q)$ bits,
which is polynomial in $\log q$. This can be improved to $O(\log^{2+\epsilon}
q)$ bits by using the algorithm of~\cite{sutherland-point-counting} to directly
compute $\Phi_\ell(j(\E),Y)$ in a space-efficient manner.

%%%%%%%%%%%%%%%%%%%%%%%%%%%%%%%%%%%%%%%%%%%%%%%%%%%%%%%%%%%%%%%%

\section{Computations}\label{sec:comput}

To test our generic low-memory algorithm for finding short
product representations in a practical setting, we implemented
black-boxes for three types of finite groups:
\begin{enumerate}
\item
$G=\E(\FF_p)$, the elliptic curve $\E:y^2=x^3+x+1$ over a finite field $\FF_p$.
\item
$G=\cl\ord$, where $\ord$ is an order in an imaginary quadratic field.\footnote{We
identify $\ord$ by its discriminant $D$ and may write $\cl{D}$ instead of $\cl\ord$.}
\item
$G=\gltwo{\FF_p}$, the group of invertible $2\times 2$ matrices over $\FF_p$.
\end{enumerate}
To simplify the implementation, we restricted to cases where $\FF_p$ is a prime field.
The groups $\E(\FF_p)$ are abelian groups, either cyclic or the product of
two cyclic groups.  The groups $\cl\ord$ are also abelian, but may be highly
non-cyclic (we specifically chose some examples with large $2$-rank), while the
groups $\gltwo{\FF_p}$ are non-abelian.

For the groups $\E(\FF_p)$, we used the
sequence of points $S=(P_1,\ldots,P_k)$ with $P_i=(x_i,y_i)$, where $x_i$ is the
$i$\textsuperscript{th} smallest positive integer for which $x_i^3+x_i+1$
is a quadratic residue $y_i^2$ modulo $p$ with $y_i \le (p-1)/2$; our target
$z$ was the point $P_{k+1}$.  For the groups $\cl\ord$, we used
the sequence $S_k$ defined in Section~\ref{sec:Sk}
with $z=[\alpha_{k+1}]$.  For the groups $\gltwo{\FF_p}$, we simply
chose a sequence $S$ of length $k$ and a target element $z$ at random.

Table~\ref{table:comput} lists performance data obtained by
applying our Pollard-$\rho$ algorithm to various groups~$G$
and sequences $S$ of densities $d=k/\log_2 n$ ranging
from just under~$2$ to slightly more than~$4$.
Each row compares expected
values with actual results that are averages over at least $10^3$ runs.

The parameter $c$ counts the number of collisions $\phi^{(i+j)}(w)=\phi^{(j)}(w)$
that were needed for a run of the algorithm to obtain a short product representation.
Typically $c$ is greater than $1$ because
not every collision yields a short product representation.  The
parameter~$\rhotot$ is the sum of $\rho=i+j$ over the $c$ collisions required,
and represents a lower bound on the number of times the map $\phi$ was evaluated. 
With efficient collision detection, the actual number is very close to~$\rhotot$
(using the method of distinguished points we were able to stay within $1\%$).

\begin{table}[p]
\begin{center}
\Small
\begin{tabular}{@{}lrrrcrrcrr@{}}
&&&&&\multicolumn{2}{c}{expected}&&\multicolumn{2}{c}{observed}\\
\cline{6-7}\cline{9-10}
$G$&$\log_2 n$& $k$&$d$&&$c$&\hspace{24pt}$\rhotot$&&$c$&\hspace{24pt}$\rhotot$\\[3pt]
\hline\\[-5pt]
$\E/\FF_{2^{20}+7}$  & 20.00& 40& 2.00&& 3.00&    3144&& 3.00&    3162\\
                           && 60& 3.00&& 2.00&    2568&& 2.01&    2581\\
                           && 80& 4.00&& 2.00&    2567&& 2.01&    2565\\[1pt]
$\E/\FF_{2^{24}+43}$ & 24.00& 48& 2.00&& 3.00&   12577&& 3.02&   12790\\
                           && 72& 3.00&& 2.00&   10269&& 2.03&   10381\\
                           && 96& 4.00&& 2.00&   10268&& 2.00&   10257\\[1pt]
$\E/\FF_{2^{28}+3}$  & 28.00& 56& 2.00&& 3.00&   50300&& 2.95&   49371\\
                           && 84& 3.00&& 2.00&   41070&& 2.02&   41837\\
                           &&112& 4.00&& 2.00&   41069&& 1.98&   40508\\[1pt]
$\E/\FF_{2^{32}+15}$ & 32.00& 64& 2.00&& 3.00&  201196&& 3.06&  205228\\
                           && 96& 3.00&& 2.00&  164276&& 1.96&  160626\\
                           &&128& 4.00&& 2.00&  164276&& 2.04&  169595\\[1pt]
$\E/\FF_{2^{36}+31}$ & 36.00& 72& 2.00&& 3.00&  804776&& 2.95&  796781\\
                           &&108& 3.00&& 2.00&  657097&& 2.00&  655846\\
                           &&144& 4.00&& 2.00&  657097&& 1.98&  657097\\[1pt]
$\E/\FF_{2^{40}+15}$ & 40.00& 80& 2.00&& 3.00& 3219106&& 2.90& 3120102\\
                           &&120& 3.00&& 2.00& 2628390&& 1.97& 2604591\\
                           &&160& 4.00&& 2.00& 2628390&& 2.06& 2682827\\[3pt]
$\cl{1-2^{40}}$      & 19.07& 40& 2.10&& 2.52&    2088&& 2.44&    2082\\
                           && 60& 3.15&& 2.00&    1859&& 2.02&    1845\\
                           && 80& 4.20&& 2.00&    1858&& 2.01&    1863\\[1pt]
$\cl{1-2^{48}}$      & 23.66& 48& 2.03&& 2.79&   10800&& 2.75&   10662\\
                           && 72& 3.04&& 2.00&    9140&& 1.97&    8938\\
                           && 96& 4.06&& 2.00&    9140&& 1.99&    9079\\[1pt]
$\cl{1-2^{56}}$      & 27.54& 56& 2.03&& 2.73&   40976&& 2.69&   40512\\
                           && 84& 3.05&& 2.00&   35076&& 2.06&   36756\\
                           &&112& 4.07&& 2.00&   35076&& 1.98&   35342\\[1pt]
$\cl{1-2^{64}}$      & 30.91& 64& 2.07&& 2.47&  125233&& 2.59&  131651\\
                           && 96& 3.11&& 2.00&  112671&& 1.98&  111706\\
                           &&128& 4.14&& 2.00&  112671&& 1.99&  111187\\[1pt]
$\cl{1-2^{72}}$      & 35.38& 72& 2.04&& 2.65&  609616&& 2.60&  598222\\
                           &&108& 3.05&& 2.00&  529634&& 2.00&  534639\\
                           &&144& 4.07&& 2.00&  529634&& 2.00&  532560\\[1pt]
$\cl{1-2^{80}}$      & 39.59& 80& 2.02&& 2.76& 2680464&& 2.80& 2793750\\
                           &&120& 3.03&& 2.00& 2283831&& 2.01& 2318165\\
                           &&160& 4.04&& 2.00& 2283831&& 2.04& 2364724\\[3pt]
$\gltwo{\FF_{37}}$   & 20.80& 42& 2.02&& 2.87&    4053&& 2.84&    4063\\
                           && 62& 2.98&& 2.00&    3384&& 1.99&    3358\\
                           && 84& 4.04&& 2.00&    3384&& 1.97&    3388\\[1pt]
$\gltwo{\FF_{67}}$   & 24.24& 48& 1.98&& 3.18&   14087&& 3.08&   13804\\
                           && 72& 2.97&& 2.00&   11168&& 2.10&   11590\\
                           && 96& 3.96&& 2.00&   11167&& 2.01&   11167\\[1pt]
$\gltwo{\FF_{131}}$  & 28.12& 56& 1.99&& 3.09&   53251&& 3.03&   52070\\
                           && 84& 2.99&& 2.00&   42851&& 1.94&   42019\\
                           &&112& 3.98&& 2.00&   42851&& 1.98&   42146\\[1pt]
$\gltwo{\FF_{257}}$  & 32.02& 64& 2.00&& 3.01&  202769&& 3.03&  204827\\
                           && 96& 3.00&& 2.00&  165237&& 2.02&  165742\\
                           &&128& 4.00&& 2.00&  165237&& 2.00&  165619\\[1pt]
$\gltwo{\FF_{511}}$  & 36.10& 72& 1.99&& 3.07&  842191&& 3.18&  886141\\
                           &&108& 2.99&& 2.00&  679748&& 1.97&  668416\\
                           &&144& 3.99&& 2.00&  679747&& 2.04&  703877\\[1pt]
$\gltwo{\FF_{1031}}$ & 40.04& 80& 2.00&& 3.03& 3276128&& 2.99& 3243562\\
                           &&120& 3.00&& 2.00& 2663155&& 2.02& 2677122\\
                           &&160& 4.00&& 2.00& 2663154&& 2.08& 2708512\\[3pt]
\end{tabular}
\end{center}
\caption{Comparison of expected vs. observed values on various groups.}
\label{table:comput}
\end{table}

The expected values of $c$ and $\rhotot$ listed in Table~\ref{table:comput} were computed
under the heuristic assumption that $\hmap:G\to\CC$ and $\pimap:\CC\to G$
are both random functions.  This implies that while iterating $\phi$
we are effectively performing simultaneous independent random walks
on $G$ and $\CC$.  Let $X$ and $Y$ be independent random
variables for the number of steps these walks 
take before reaching a collision, respectively. 
The probability that $\pi(s)=\pi(t)$ in Step~5 is $P(X \le Y)$,
and the algorithm then proceeds to find a short product representation with
probability $1/2$.

Using the probability density $u\exp(-u^2/2)du$ of
$X/\sqrt{\#G}$ and $Y/\sqrt{\#\CC}$, we find
\[
\Exp{c} = 2/{P(X\le Y)} = 2(1+r),
\]
where $r=\#G/\#\CC$.  One may also compute
\[
\Exp{\rhotot} = \Exp{c}\Exp{\min(X,Y)} = \sqrt{2\pi n(1+r)}.
\]
For $d > 2$, we have $r\approx 0$ for large $n$,
so that $\Exp{c}\approx 2$ and $\Exp{\rhotot}\approx \sqrt{2\pi n}$.
For $d=2$, we have $\Exp{c}=3$ and $\Exp{\rhotot}=\sqrt{3\pi n}$
(when $k$ is even).  For $d < 2$, the value of $\Exp{c}$ increases with $n$
and we have $\Exp{\rhotot}=O(n^{(4-d)/4})$.

\medskip

In addition to the tests summarized in Table~\ref{table:comput}, we applied our
low memory algorithm to some larger problems that would be
quite difficult to address with the baby-step giant-step
method.  Our first large test used $G=\E(\FF_p)$ with
$p=2^{80}+13$, which is a cyclic group of
order $n=p+1+1475321552477$, and the sequence $S=(P_1,\ldots,P_{k})$ 
with points $P_i$ defined as above with $k=200$, which gives $d\approx 2.5$.
Our target element was $z=P_{201}$ with $x$-coordinate $391$.
The computation was run in parallel on $32$ cores (3.0~GHz AMD Phenom~II), using the
distinguished points method.\footnote{In this parallel setting we may have collisions between two distinct
walks (a $\lambda$-collision), or a single walk may collide with itself (a $\rho$-collision).  Both
types are useful.}  The second collision yielded a
short product representation after evaluating the map~$\phi$
a total of $1480862431620 \approx 1.35\sqrt{n}$ times.

After precomputing $655360$ partial products (as discussed in Section~\ref{sec:analysis}),
each evaluation of $\phi$ used $5$ group operations, compared to an average of $50$
without precomputation, and this required just $10$ megabytes of memory.
The entire computation used approximately $140$~days of CPU time, 
and the elapsed time was about $4$~days.
We obtained a short product representation for $z$ as the sum of $67$ points $P_i$ with $x$-coordinates less than $391$.
In hexadecimal notation, the bit-string that identifies the corresponding subsequence of $S$ is:
\begin{center}
\texttt{542ab7d1f505bdaccdbeb6c2e92180d5f38a20493d60f031c1}
\end{center}

\medskip

Our second large test used the group $G=\cl{1-2^{160}}$, which is
isomorphic to
\[
(\ZZ/2\ZZ)^{8} \times \ZZ/4\ZZ  \times \ZZ/8\ZZ \times \ZZ/80894875660895214584\ZZ,
\]
see~\cite[Table B.4]{sutherland:thesis}.
We used the sequence $S_k$ with $k=200$, and
chose the target $z=[\alpha_{201}]$ with $\N(\alpha_{201})=2671$.
We ran the computation in parallel on $48$ cores, and needed
$3$ collisions to obtain a short product representation, which
involved a total of $2856153808020\approx 3.51\sqrt{n}$ evaluations of $\phi$.
As in the first test, we precomputed $655360$ partial products so that
each evaluation of $\phi$ used $5$ group operations.  Approximately $900$ days
of CPU time were used (the group operation in $\cl{D}$ is slower
than in the group $E(\FF_p)$ used in our first example).  We obtained
a representative for the ideal class $z$ as the product of $106$
ideals with prime norms less than $2671$.
The bit-string that encodes the corresponding subsequence of $S_k$ is:
\begin{center}
\texttt{5cf854598d6059f607c6f17b8fb56314e87314bee7df9164cd}
\end{center}

%%%%%%%%%%%%%%%%%%%%%%%%%%%%%%%%%%%%%%%%%%%%%%%%%%%%%%%%%%%%%%%%

\section*{Acknowledgments}

The authors are indebted to Andrew Shallue for his kind help and advice in
putting our result in the context of subset sum problems, and to Steven Galbraith
for his useful feedback on an early draft of this paper.

%%%%%%%%%%%%%%%%%%%%%%%%%%%%%%%%%%%%%%%%%%%%%%%%%%%%%%%%%%%%%%%%

\bibliographystyle{plain}

\begin{thebibliography}{10}

\bibitem{alon-barak-manber}
Noga Alon, Amnon Barak, and Udi Manber.
\newblock On disseminating information reliably without broadcasting.
\newblock In Radu Popescu-Zeletin, Gerard {Le Lann}, and Kane~H. Kim, editors,
  {\em Proceedings of the 7\textsuperscript{th} International Conference on
  Distributed Computing Systems}, pages 74--81. IEEE Computer Society Press,
  1987.

\bibitem{alon-milman}
Noga Alon and Vitali~D. Milman.
\newblock {$\lambda_1$}, isoperimetric inequalities for graphs, and
  superconcentrators.
\newblock {\em Journal of Combinatorial Theory, Series B}, 38:73--88, 1985.

\bibitem{babai-erdos}
L\'{a}szl\'{o} Babai and Paul Erd\H{o}s.
\newblock Representation of group elements as short products.
\newblock {\em North-Holland Mathematics Studies}, 60:27--30, 1982.

\bibitem{bach-erh}
Eric Bach.
\newblock Explicit bounds for primality testing and related problems.
\newblock {\em Mathematics of Computation}, 55(191):355--380, 1990.

\bibitem{bisson-grh}
Gaetan Bisson.
\newblock Computing endomorphism rings of elliptic curves under the {GRH},
  2010.
\newblock In preparation.

\bibitem{bisson-sutherland}
Gaetan Bisson and Andrew~V. Sutherland.
\newblock Computing the endomorphism ring of an ordinary elliptic curve over a
  finite field.
\newblock {\em Journal of Number Theory}, Special Issue on Elliptic Curve
  Cryptography, 2009.
\newblock To appear.

\bibitem{brent}
Richard~P. Brent.
\newblock An improved {Monte Carlo} factorization algorithm.
\newblock {\em BIT Numerical Mathematics}, 20:176--184, 1980.

\bibitem{quantum-iso}
Andrew~M. Childs, David Jao, and Vladimir Soukharev.
\newblock Constructing elliptic curve isogenies in quantum subexponential time,
  2010.
\newblock Preprint available at \url{http://arxiv.org/abs/1012.4019}.

\bibitem{eggleton-erdos}
Roger~B. Eggleton and Paul Erd\H{o}s.
\newblock Two combinatorial problems in group theory.
\newblock {\em Acta Arithmetica}, 28:247--254, 1975.

\bibitem{erdos-renyi}
Paul Erd\H{o}s and Alfr\'{e}d R\'{e}nyi.
\newblock Probabilistic methods in group theory.
\newblock {\em Journal d'Analyse Math\'{e}matique}, 14(1):127--138, 1965.

\bibitem{galbraith}
Steven~D. Galbraith.
\newblock Constructing isogenies between elliptic curves over finite fields.
\newblock {\em Journal of Computational Mathematics}, 2:118--138, 1999.

\bibitem{galbraith-hess-smart}
Steven~D. Galbraith, Florian Hess, and Nigel~P. Smart.
\newblock Extending the {GHS} {Weil} descent attack.
\newblock In Lars~R. Knudsen, editor, {\em Advances in Cryptology--EUROCRYPT
  '02}, volume 2332 of {\em Lecture Notes in Computer Science}, pages 29--44.
  Springer, 2002.

\bibitem{hafner-mccurley}
James~L. Hafner and Kevin~S. McCurley.
\newblock A rigorous subexponential algorithm for computing in class groups.
\newblock {\em Journal of the American Mathematical Society}, 2(4):837--850,
  1989.

\bibitem{howgravegraham-joux}
Nick Howgrave-Graham and Antoine Joux.
\newblock New generic algorithms for hard knapsacks.
\newblock In Henri Gilbert, editor, {\em Advances in Cryptology--EUROCRYPT
  '10}, volume 6110 of {\em Lecture Notes in Computer Science}, pages 235--256.
  Springer, 2010.

\bibitem{impagliazzo-naor}
Russel Impagliazzo and Moni Naor.
\newblock Efficient cryptographic schemes provably as secure as subset sum.
\newblock {\em Journal of Cryptology}, 9(4):199--216, 1996.

\bibitem{expander-grh}
David Jao, Stephen~D. Miller, and Ramarathnam Venkatesan.
\newblock Expander graphs based on {GRH} with an application to elliptic curve
  cryptography.
\newblock {\em Journal of Number Theory}, 129(6):1491--1504, 2009.

\bibitem{karp}
Richard~M. Karp.
\newblock Reducibility among combinatorial problems.
\newblock In Raymond~E. Miller, James~W. Thatcher, and Jean~D. Bohlinger,
  editors, {\em Complexity of Computer Computations}, pages 85--103. Plenum
  Press, 1972.

\bibitem{knuth-art2}
Donald~E. Knuth.
\newblock {\em The Art of Computer Programming, Volume II: Seminumerical
  Algorithms}.
\newblock Addison-Wesley, 1998.

\bibitem{knuth-art4f2}
Donald~E. Knuth.
\newblock {\em The Art of Computer Programming, Volume IV, Fascicle 2:
  Generating all Tuples and Permutations}.
\newblock Addison-Wesley, 2005.

\bibitem{merkle-hellman}
Ralph Merkle and Martin Hellman.
\newblock Hiding information and signatures in trapdoor knapsacks.
\newblock {\em IEEE Transactions on Information Theory}, 24(5):525--530, 1978.

\bibitem{pollard}
John~M. Pollard.
\newblock A {Monte Carlo} method for factorization.
\newblock {\em BIT Numerical Mathematics}, 15(3):331--334, 1975.

\bibitem{schonhage-fastforms}
Arnold Sch{\"o}nhage.
\newblock Fast reduction and composition of binary quadratic forms.
\newblock In Stephen~M. Watt, editor, {\em International Symposium on Symbolic
  and Algebraic Computation--ISSAC '91}, pages 128--133. ACM Press, 1991.

\bibitem{schoof-pointcounting}
Ren{\'e} Schoof.
\newblock Counting points on elliptic curves over finite fields.
\newblock {\em Journal de Th{\'e}orie des Nombres de Bordeaux}, 7:219--254,
  1995.

\bibitem{schroeppel-shamir}
Richard Schroeppel and Adi Shamir.
\newblock A {$T=O(2^{n/2}), S=O(2^{n/4})$} algorithm for certain {NP}-complete
  problems.
\newblock {\em SIAM Journal of Computing}, 10(3):456--464, 1981.

\bibitem{sedgewick}
Robert Sedgewick and Thomas~G. Szymanski.
\newblock The complexity of finding periods.
\newblock In {\em Proceedings of the 11\textsuperscript{th} ACM Symposium on
  the Theory of Computing}, pages 74--80. ACM Press, 1979.

\bibitem{shoup}
Victor Shoup.
\newblock Lower bounds for discrete logarithms and related problems.
\newblock In {\em Advances in Cryptology--EUROCRYPT '97}, volume 1233 of {\em
  Lecture Notes in Computer Science}, pages 256--266. Springer-Verlag, 1997.
\newblock Revised version.

\bibitem{siegel}
Carl~Ludwig Siegel.
\newblock {\"U}ber die {C}lassenzahl quadratischer {Z}ahlk{\"o}rper.
\newblock {\em Acta Arithmetica}, 1:83--86, 1935.

\bibitem{sobol-random-sequences}
Ilya~M. Sobol.
\newblock On periods of pseudo-random sequences.
\newblock {\em Theory of Probability and its Applications}, 9:333--338, 1964.

\bibitem{sutherland-point-counting}
Andrew~V. Sutherland.
\newblock Genus 1 point counting in quadratic space and essentially quartic
  time.
\newblock in preparation.

\bibitem{sutherland:thesis}
Andrew~V. Sutherland.
\newblock Order computations in generic groups.
\newblock {PhD} thesis, MIT, 2007.
\newblock \url{http://groups.csail.mit.edu/cis/theses/sutherland-phd.pdf}.

\bibitem{teske}
Edlyn Teske.
\newblock A space efficient algorithm for group structure computation.
\newblock {\em Mathematics of Computation}, 67:1637--1663, 1998.

\bibitem{vanoorschot-wiener}
Paul~C. van Oorschot and Michael~J. Wiener.
\newblock Parallel collision search with cryptanalytic applications.
\newblock {\em Journal of Cryptology}, 12:1--28, 1999.

\bibitem{white}
Edward White.
\newblock Ordered sums of group elements.
\newblock {\em Journal of Combinatorial Theory, Series A}, 24:118--121, 1978.

\end{thebibliography}

\end{document}